\documentclass[preprint,3p,12pt]{elsarticle}

\usepackage{lineno,hyperref}
\modulolinenumbers[5]

\usepackage{hyperref}

\usepackage{epstopdf}

\usepackage{mathrsfs}
\usepackage{amssymb}
\usepackage{amsfonts}
\usepackage{latexsym}
\usepackage{amsmath}
\usepackage{xcolor}
\usepackage{wrapfig}
\usepackage{floatflt}

\usepackage{graphicx}
\def\lb{\label}

\newcommand{\er}[1]{\textrm{(\ref{#1})}}

\newtheorem{theorem}{\bf Theorem}[section]
\newtheorem{lemma}[theorem]{\bf Lemma}
\newtheorem{definition}[theorem]{\bf Definition}

\def\a{\alpha}  \def\cA{{\mathcal A}}       
\def\b{\beta}   \def\cB{{\mathcal B}}       \def\mB{{\mathscr B}}
\def\g{\gamma}  \def\cC{{\mathcal C}}       \def\mC{{\mathscr C}}
  \def\cD{{\mathcal D}}       
\def\d{\delta}         
         
    \def\cG{{\mathcal G}}       \def\mG{{\mathscr G}}
          \def\mH{{\mathscr H}}
    \def\cI{{\mathcal I}}

         \def\mL{{\mathscr L}}
\def\l{\lambda}        
        
\def\m{\mu}

\def\s{\sigma}         
         
\def\t{\tau}

\def\ve{\varepsilon}           

\def\Z{{\mathbb Z}}       \def\C{{\mathbb C}}    
    \def\N{{\mathbb N}}   

\def\ol{\overline}               \def\wt{\widetilde}

\let\ge\geqslant                 \let\le\leqslant

\def\iy{\infty}
\def\sm{\setminus}               \def\es{\emptyset}
\def\ss{\subset}                 \def\ts{\times}

\def\el2{\ell^{\,2}}             \def\1{1\!\!1}

\def\det{\mathop{\mathrm{det}}\nolimits}

\def\BBox{\hspace{1mm}\vrule height6pt width5.5pt depth0pt \hspace{6pt}}

\let\ge\geqslant
\let\le\leqslant

\newcommand{\ca}{\begin{cases}}
\newcommand{\ac}{\end{cases}}
\newcommand{\ma}{\begin{pmatrix}}
\newcommand{\am}{\end{pmatrix}}
\def\eq{\begin{equation}}
\def\qe{\end{equation}}
\def\[{\begin{equation}}
\def\]{\end{equation}}

\def\BBox{\hspace{1mm}\vrule height6pt width5.5pt depth0pt \hspace{6pt}}

\begin{document}

\begin{frontmatter}

\title{On the extension of Fredholm determinants to the mixed multidimensional
integral operators with regulated kernels}
\date{\today}

\author
{Anton A. Kutsenko}


\address{Jacobs University (International University Bremen), 28759 Bremen, Germany; email: akucenko@gmail.com}

\begin{abstract}
We extend the classical trace (and determinant) known for the
integral operators
$$
 ({\mathcal I}+)\int_{[0,1)^N}{\bf A}({\bf k},{\bf
 x}){\bf u}({\bf x})d{\bf x}
$$
with matrix-valued kernels ${\bf A}$
to the operators of the form
$$
 \sum_{\alpha}\int_{[0,1)^{|\alpha|}}{\bf A}({\bf
 k},{\bf x}_{\alpha}){\bf u}({\bf k}_{\overline{\alpha}},{\bf x}_{\alpha})d{\bf
 x}_{\alpha},
$$
where $\alpha$ are arbitrary subsets of the set $\{1,...,N\}$. Such
operators form a Banach algebra containing simultaneously all
integral operators of the dimensions $\leqslant N$. In this sense,
it is a largest algebra where explicit traces and determinants are
constructed. Such operators arise naturally in the mechanics and
physics of waves propagating through periodic structures with
various defects. We give an explicit representation of the inverse
operators (resolvent) and describe the spectrum by using zeroes of
the determinants. Due to the structure of the operators, we have
$2^N$ different determinants, each of them describes the spectral
component of the corresponding dimension.
\end{abstract}



\end{frontmatter}


\section{Introduction}

Defects in periodic structures and corresponding periodic operators
play a significant role in various fields of science, see, e.g.
\cite{WHMH,CWC,JMFVJH,JJ,K3}. It is shown in \cite{Kjmaa,Kjmaa1}
that periodic operators with crossing defects (periodic sublattices)
of various dimensions are unitarily equivalent to the integral
operators of the special form. We start with introducing these
integral operators. At first, we consider the operators with
staircase (piecewise constant) kernels because such kernels
significantly simplify the proof of key Lemma \ref{L7} and they are
important in numerical applications. Finally, we generalize our
results to the continuous (or even regulated) kernels. Let $N,M$ be
some positive integer numbers. Introduce the Hilbert space
$L^2_{N,M}$ of square-integrable vector-valued (if $M>1$) functions
acting on the cube $[0,1)^N$ and integrals $\langle...\rangle_{\a}$
by
\[\lb{001}
 L^2_{N,M}:=L^2([0,1)^N,\C^M),\ \
 \langle...\rangle_{\a}:=\int_{[0,1)^{|\a|}}...
 d{\bf x}_{\a}
\]
where $\a\in I_N=2^{\{1,..,N\}}$ is a subset of the set
$\{1,..,N\}$, $|\a|$ is the number of elements of $\a$, and $d{\bf
x}_{\a}=\prod_{i\in\a}dx_i$, and $\langle{\bf a}\rangle_{\es}:={\bf
a}$ for any ${\bf a}$. The dots $...$ in \er{001} means any
matrix-valued function depending on ${\bf
x}_{\a}=(x_i)_{i\in\a}\in[0,1)^{|\a|}$, the components $x_i$ in
${\bf x}_{\a}$ are arranged in the order of increasing indices. Let
$\a$, $\b$ be disjoint subsets from $I_N$. It is convenient to use
the following notation
\[\lb{002}
 {\bf y}_{\b}\diamond{\bf x}_{\a}={\bf x}_{\a}\diamond{\bf y}_{\b}={\bf z}_{\a\cup\b},\ \ where\ \
 z_i=\ca x_i,& i\in\a,\\ y_i,& i\in\b.\ac
\]
All components of ${\bf x}_{\a},{\bf y}_{\b},{\bf z}_{\a\cup\b}$ in
\er{002} are arranged in the order of increasing indices. The
operation $\diamond$ is also associative. Define the following
subset of the algebra $\mB_{N,M}$ of all bounded operators acting on
$L^2_{N,M}$. Let $h=1/p$ with some $p\in\N$. Define $\chi^h_i(y)$ is
a function which is $1$ for $y\in[i/p,(i+1)/p)$ and is $0$
otherwise. We call the function $A(y_1,...,y_R)$ as $h$-staircase if
it is a linear combination of $\prod_{j=1}^R\chi^h_{i_j}(y_j)$. We
call the matrix-valued function ${\bf A}(y_1,...,y_M)$ as
$h$-staircase if all entries are $h$-staircase.

\begin{definition} \lb{D1}
Let $\mL_{N,M}^h$ be the set
\[\lb{003}
 \mL_{N,M}^h=\{\cA:\ \cA=\sum_{\a\in I_N}\langle{\bf A}_{\a}({\bf k},{\bf
 x}_{\a})\cdot\rangle_{\a}\},
\]
where ${\bf A}_{\a}$ are any $h$-staircase $M\ts M$ matrix-valued
functions depending on ${\bf k}=(k_i)_{i=1}^N\in[0,1)^N$, ${\bf
x}_{\a}=(x_i)_{i\in\a}\in[0,1)^{|\a|}$. The dot $\cdot$ denotes the
place of the operator argument, i.e.
\[\lb{004a}
 \cA{\bf u}({\bf k})=\sum_{\a\in I_N}\langle{\bf A}_{\a}({\bf k},{\bf
 x}_{\a}){\bf u}({\bf k}_{\ol{\a}}\diamond{\bf x}_{\a})\rangle_{\a}, \ \
{\bf u}\in L^2_{N,M},
\]
where $\ol{\a}$ is the complement to the set $\a\in I_N$.
\end{definition}

It is obvious that $\mL$ (sometimes we will omit indices $N,M$ and
$h$ for convenience) is a linear subspace of $\mB$. Let us introduce
the following positive function on $\mL$.

\begin{definition}\lb{D2}
Denote
\[\lb{004}
 \|\cA\|_{\mL}=\sum_{\a}\max_{({\bf k},{\bf
 x}_{\a}),i}\sum_{j=1}^M|a_{i,j,\a}({\bf k},{\bf
 x}_{\a})|,
\]
where ${\bf A}_{\a}=(a_{i,j,\a})_{i,j=1}^M$. Due to Lemma \ref{L1},
definition \er{004} is correct.
\end{definition}

The next proposition describes the basic properties of $\mL$.

\begin{theorem}\lb{T1} The function $\|...\|_{\mL}$ is a norm on $\mL$. This norm is stronger than the standard
operator norm $\|...\|_{\mB}$. The convergence in $\|...\|_{\mL}$ is
equivalent to the uniform convergence of the coefficients ${\bf
A}_{\a}$. The structure $(\mL,\|...\|_{\mL})$ is a Banach algebra
(associative, non-commutative).
\end{theorem}

While the sum of two operators from $\mL$ leads to the sum of the
corresponding matrix-valued functions in their representations, the
product (composition) is more complicated. For example, it is not
difficult to check the following identity for the product of two
terms
\[\lb{no1}
 \langle{\bf A}({\bf k},{\bf x}_{\a})\cdot\rangle_{\a}\langle{\bf B}({\bf k},{\bf
 x}_{\b})\cdot\rangle_{\b}=\langle({\bf A}\circ{\bf B})({\bf k},{\bf
 x}_{\a\cup\b})\cdot\rangle_{\a\cup\b},
\]
where the matrix-valued function ${\bf A}\circ{\bf B}$ is defined by
\[\lb{no2}
 ({\bf A}\circ{\bf B})({\bf k},{\bf x}_{\a\cup\b})=\int_{[0,1)^{|\a\cap\b|}}{\bf
 A}({\bf k},{\bf x}_{\a\sm\b}\diamond{\bf z}_{\a\cap\b})
 {\bf B}({\bf k}_{\ol{\a}}\diamond{\bf z}_{\a\cap\b}\diamond{\bf x}_{\a\sm\b},{\bf
 x}_{\b})d{\bf z}_{\a\cap\b}.
\]
Now, define the "building blocks" of $\mL$ and ${\rm Inv}(\mL)$
(invertible operators). We assume the invertibility in $\mB$ but we
show below that if $\cA\in\mL$ is invertible in $\mB$ then
$\cA^{-1}\in\mL$. This is more or less evident for
finite-dimensional $\mL$. Introduce
\[\lb{012}
 \mL_{\a}=\{\langle{\bf A}({\bf k},{\bf x}_{\a})\cdot\rangle_{\a}:\ \forall{\bf
 A}\}\ss\mL;\ \ \mG_{\es}={\rm Inv}(\mL_{\es}),\ \ \mG_{\a}={\rm
 Inv}(\cI+\mL_{\a})\ if\ \a\ne\es,
\]
where $\cI:={\bf I}\cdot$ is the identity operator and ${\bf I}$ is
the identity matrix.
\begin{theorem}\lb{T3} The sets $\mL_{\a}$ are closed subalgebras of $\mL$, the sets
$\mG_{\a}$ are closed subgroups of ${\rm Inv}(\mL)$, satisfying
\[\lb{013a}
 \mL_{\a}\cap\mL_{\b}=\{0\},\ \ \mG_{\a}\cap\mG_{\b}=\{\cI\}\ \ for
 \ \ \a\ne\b.
\]
The following identities hold true
\[\lb{013}
 \mL_{\a}\mL_{\b}\ss\mL_{\a\cup\b},\ \ \ \mL=\sum_{\a\in I_N}\mL_{\a},\ \ \ {\rm Inv}(\mL)=\prod_{\a\in
 I_N}\mG_{\a},
\]
where the order of terms in the sum is not important but the terms
in the product are arranged in ascending (or descending) order of
$|\a|$. Moreover, for given $\cA$ and for given order of terms the
corresponding representations as the sum and the product of
elementary operators are unique.
\end{theorem}

The explicit procedure of finding the components in the product
\er{013} is given in the proof of Theorem \ref{T3}. Because the
inverse of the elementary operators from the product have an
explicit form, we obtain the explicit form for the inverse of the
product of operators. This means that having $\cA$ we can check its
invertibility and find $\cA^{-1}$ explicitly. Before introducing the
trace and the determinant, let us introduce the commutative algebras
of complex scalar $h$-staircase functions acting on different cubes
$[0,1)^r$
\[\lb{007}
 \mC^h_{\a,N}=\{f({\bf k}_{\ol{\a}}),{\bf
 k}_{\ol{\a}}=(k_i)_{i\not\in\a}\in[0,1)^{N-|\a|}\},
\ \
 \mC^h_N=\bigoplus_{\a\in
 I_N}\mC^h_{\a,N}.
\]
The algebra $\mC_{\es,N}$ is the algebra of all complex
$h$-staircase scalar functions defined on the cube $[0,1)^N$. If
$\a=\{1,..,N\}$ then we put $\mC_{\a,N}:=\C$. All $\mC_{\a,N}$ are
subalgebras of $\mC_{\es,N}$ and they consist of the functions
independent on some variables $k_i$. In this sense if $|\a|=|\b|$
but $\a\ne\b$ then $\mC_{\a}\ne\mC_{\b}$. All operations $+$, $-$,
$*$, $\ol{\bf f}$ (complex conjugation), $\exp({\bf f})$, etc in the
algebra $\mC$ are assumed to be componentwise. Taking the standard
norm $\|{\bf f}\|=\max_{\a,{\bf k}_{\ol{\a}}}|f_{\a}({\bf
k}_{\ol{\a}})|$ we construct the commutative finite-dimensional
Banach algebra $(\mC,\|...\|_{\mC})$. Introduce the following
mappings.
\begin{definition}\lb{D3}
Define the following mapping
\[\lb{008}
 \pmb{\t}:\mL\to\mC,\ \ \pmb{\t}(\cA)=(\t_{\a}(\cA)),\ \
 \t_{\a}(\cA)=\t_{\a}(\cA)({\bf k}_{\ol{\a}}):={\rm Tr}\langle{\bf A}_{\a}({\bf k}_{\ol{\a}}\diamond{\bf x}_{\a},{\bf
 x}_{\a})\rangle_{\a}.
\]
Now, we fix some ascending (or descending) order of $|\a|$. Below we
show that the definition does not depend on the order. If
$\cA\in{\rm Inv}(\mL)$ then by Theorem \ref{T3} $\cA=\prod_{\a\in
I_N}\cG_{\a}$, where $\cG_{\a}=\cI+\langle{\bf G}_{\a}({\bf k},{\bf
x}_{\a})\cdot\rangle_{\a}$. We can write $\cG_{\a}=\cG_{\a}({\bf
k}_{\ol{\a}})$, where for any fixed ${\bf k}_{\ol{\a}}$ the operator
$\cG_{\a}({\bf k}_{\ol{\a}})=\cI+\langle{\bf G}_{\a}({\bf
k}_{\ol{\a}}\diamond{\bf k}_{\a},{\bf x}_{\a})\cdot\rangle_{\a}$ (or
$\cG_{\es}({\bf k})={\bf G}_{\es}({\bf k})\cdot$ for $\a=\es$) is a
finite-rank operator acting on $L^2([0,1)^{|\a|},\C^M)$. Hence, we
can define the determinant $\pi_{\a}$ of such operators
\[\lb{008b}
\pi_{\a}(\cG_{\a})=\pi_{\a}(\cG_{\a})({\bf
k}_{\ol{\a}}):=\pi_{\a}(\cG_{\a}({\bf k}_{\ol{\a}}))
\]
This leads to the mapping
\[\lb{008a}
 \pmb{\pi}:{\rm Inv}(\mL)\to{\rm Inv}(\mC),\ \ \pmb{\pi}(\cA)=(\pi_{\a}(\cG_{\a})),
\]
since $\cG_{\a}({\bf k}_{\ol{\a}})$, and hence their determinants,
are staircase functions. The set ${\rm Inv}(\mC)$ consists of all
functions from $\mC$ that have no zeroes, since the determinant of
finite-rank invertible operators is not zero.
\end{definition}

{\bf Remark (computation of $\pmb{\pi}$).} While $\pmb{\t}$ is
simple, the computation of $\pmb{\pi}$ is more complicated.
Nevertheless, there are various formulas for $\pmb{\pi}$. Separating
the variables (see the proof of Lemma \ref{L2}) $\cI+\langle{\bf
G}({\bf k},{\bf x}_{\a})\cdot\rangle_{\a}=\cI+{\bf C}({\bf
k})\langle{\bf D}({\bf x}_{\a})\cdot\rangle_{\a}$ we can compute the
determinant
\[\lb{C1}
 \pi_{\a}(\cI+\langle{\bf G}({\bf k},{\bf
 x}_{\a})\cdot\rangle_{\a})=\det{\bf E}({\bf k}_{\ol{\a}}),\ \ {\bf E}({\bf k}_{\ol\a})={\bf I}+\langle{\bf D}({\bf
 x}_{\a}){\bf C}({\bf k}_{\ol{\a}}\diamond{\bf x}_{\a})\rangle_{\a}.
\]
This is a very convenient computation of the determinant of
finite-rank operators, see, e.g. \cite{GGK}, \cite{K1a}. We can use
the Fredholm formulas
\[\lb{F1}
 \pi_{\a}(\cI+\langle{\bf G}_{\a}({\bf k},{\bf
x}_{\a})\cdot\rangle_{\a})=1+\sum_{n=1}^{+\iy}\frac1{n!}\int_{[0,1)^{n|\a|}}P_n({\bf
k}_{\ol{\a}},{\bf x}_{1\a},...,{\bf x}_{n\a})d{\bf x}_{1\a}...d{\bf
x}_{n\a},
\]
\[\lb{F2}
P_n({\bf k}_{\ol{\a}},{\bf x}_{1\a},...,{\bf x}_{n\a})=\det\ma{\bf
G}_{\a}({\bf k}_{\ol\a}\diamond{\bf x}_{1\a},{\bf x}_{1\a}) & ... &
{\bf G}_{\a}({\bf k}_{\ol\a}\diamond{\bf x}_{1\a},{\bf x}_{n\a}) \\
. & ... & . \\  {\bf G}_{\a}({\bf k}_{\ol\a}\diamond{\bf
x}_{n\a},{\bf x}_{1\a}) & ... & {\bf G}_{\a}({\bf
k}_{\ol\a}\diamond{\bf x}_{n\a},{\bf x}_{n\a})\am
\]
which also lead to the determinant, see, e.g., \cite{Fredholm,GGK}.

The next theorem shows us that $\pmb{\t}$ and $\pmb{\pi}$ are the
trace and the determinant (totally, i.e. not only for the elementary
operators from $\mG_{\a}$). The theory of traces and determinants
has own special interest, see, e.g. \cite{GGK,GKZ,GGRW,S,S1}. In our
case, we extend the traces and determinants to the algebra $\mL$
containing simultaneously operators of multiplication by
matrix-valued functions and various classes of integral operators.
That is why, in our case, the trace and the determinant are vectors
consisting of some functions (usually, traces and determinants are
just numbers). At the end of this section we discuss how to extend
the trace and the determinant to the operators with "general"
kernels.

\begin{theorem}\lb{T2}
For any $\l,\m\in\C$, $\cA,\cB\in\mL$, and $\cC,\cD\in {\rm
Inv}(\mL)$ the following identities hold true
\[\lb{011}
 \pmb{\t}(\l\cA+\m\cB)=\l\pmb{\t}(\cA)+\m\pmb{\t}(\cB),\ \
 \pmb{\t}(\cA\cB)=\pmb{\t}(\cB\cA),\ \
 \pmb{\pi}(\cC\cD)=\pmb{\pi}(\cC)\pmb{\pi}(\cD).
\]
Moreover $\pmb{\t}$, $\pmb{\pi}$ are continuous mappings satisfying
$\pmb{\pi}(e^{\cA})=e^{\pmb{\t}(\cA)}$ and
$\|\pmb{\t}\|_{\mL\to\mC}=M$.
\end{theorem}

Results of Theorem \ref{T2} mean that $\pmb{\t}$ and $\pmb{\pi}$ are
the trace and the determinant. They are in a good agreement with the
standard trace and determinant of finite-rank operators. In this
sense they are unique. Of course, taking linear combinations of
$\t_{\a}$ and products of $\pi_{\a}$ we can construct other traces
and determinants but only $\pmb{\pi}$ and $\pmb{\t}$ contain the
most complete information about the operator. Since
$\pmb{\pi}(e^{\cA})=e^{\pmb{\t}(\cA)}$ we can use the analog of
Plemelj-Smithies formula
\[\lb{PS}
 {\pi}_{\a}(\cI+\cA)=1+\sum_{n=1}^{+\iy}\frac{\det}{n!}
 \ma\t_{\a}(\cA) & n-1 & 0 & 0 & ... & 0 & 0 \\
                                                              \t_{\a}(\cA^2) & \t_{\a}(\cA) & n-2 & 0 & ... & 0 & 0 \\
                                                              . & . & . & . & ... & . & . \\
                                                              \t_{\a}(\cA^{n-1}) & \t_{\a}(\cA^{n-2}) & \t_{\a}(\cA^{n-3}) & . & ... & \t_{\a}(\cA) & 1 \\
                                                              \t_{\a}(\cA^{n}) & \t_{\a}(\cA^{n-1}) & \t_{\a}(\cA^{n-2}) & . & ... & \t_{\a}(\cA^2) &
                                                              \t_{\a}(\cA)\am.
\]
It is seen that the spectrum
\[\lb{023}
 \s(\cA)=\{\l\in\C:\ \not\exists(\l\cI-\cA)^{-1}\}
\]
can be determined as zeroes of the continuously extended determinant
\[\lb{024}
 \s(\cA)=\bigcup_{\a\in I_N}\{\l\in\C:\ \pi_{\a}(\l\cI-\cA)({\bf
 k}_{\ol{\a}})=0\ for\ some\ {\bf k}_{\ol{\a}}\}
\]
because if all $\pi_{\a}(\cA)$ are non-zero for some $\cA$ then
$\cA$ is invertible (see \er{L2002} and below). The identities
$\pi_{\a}(\l\cI-\cA)=0$ lead to implicit functions
$\l=\l_{\ol{\a}}({\bf k}_{\ol{\a}})$ which are very important in the
study of mechanics and physics of waves (about spectral problems in
wave dynamics see, e.g., \cite{B,Kuchment,K1}). The functions
$\l_{\ol{\a}}$ generalize the well-known Floquet-Bloch (F-B)
dispersion branches and describe the dependence of the spectral
parameter $\l$ (energy, frequency) on the wave number (or
quasimomentum) ${\bf k}_{\ol{\a}}$ corresponding to the defect
modes. Defect modes are non-attenuated (quasiperiodic $\sim e^{i{\bf
k}_{\ol{\a}}\cdot{\bf n}_{\ol{\a}}}$, where ${\bf
n}_{\ol{\a}}\in\Z^{|\ol{\a}|}$ parameterize the defect sublattice)
along the defect and exponentially decaying in the perpendicular
directions (so-called guided or Rayleigh modes). Due to the
attenuation, the dispersion branches $\l_{\ol{\a}}$ do not depend on
some components of the full quasimomentum ${\bf k}$. This is a
significant difference between the periodic structures with defects
of lower dimensions and purely periodic structures without defects,
where there are no attenuation of the modes and F-B dispersion
branches depend on the whole ${\bf k}$, i.e. we have only
$\l_{\es}({\bf k})$.

{\bf Important remark.} Note that
\[\lb{ir1}
 \mL^{\frac1n}_{N,M}\cup\mL^{\frac1m}_{N,M}\ss\mL^{\frac1{nm}}_{N,M},\
 \ n,m\in\N.
\]
So we can consider the algebra
\[\lb{ir2}
 \mL^0_{N,M}=\bigcup_{n=1}^{\iy}\mL_{N,M}^{\frac{1}{n}}.
\]
Taking the completeness of $\mL^0$ by the norm $\|\cdot\|_{\mL}$
(which is the same for all $\mL_{N,M}^{\frac{1}{n}}$) we obtain the
Banach algebra $\ol{\mL}^0_{N,M}$ of multidimensional integral
operators with regulated kernels (possible discontinuities are in
rational points). In particular, this algebra contains the operators
with continuous kernels. The same procedure leads to $\ol{\mC}^0$
consisting of regulated vector-valued functions. It is seen that the
trace $\pmb{\t}$ is continuous in $\ol{\mL}^0_{N,M}$. Hence, we can
continuously extend the determinant $\pmb{\pi}$ to the operators
from ${\rm Inv}(\ol{\mL}^0_{N,M})$. This reminds the Schmidt idea
(see, e.g., \cite{Zemyan}, 2.4 The Fredholm Theorems, p. 48) of
approximation of the integral kernel by a separable kernel and a
small kernel which leads to two integral equations with explicit
solutions. Note that taking the other systems of intervals (not only
$[i/p,i/p+1/p)$) we can obtain different classes of regulated
kernels for which the traces and determinants are defined.

\section{Proof of the results}

\begin{lemma}\lb{L1} If $\cA=0$ then all ${\bf A}_{\a}=0$, i.e. the
representation \er{002} is unique.
\end{lemma}
{\it Proof.} The proof repeats the arguments from
\cite{Kjmaa,Kjmaa1,K1a}. We give it in a short form. The proof
consists of $2^N$ steps. On the first step we suppose that ${\bf
A}_{\es}\ne0$. Then there is ${\bf k}^0=(k_1^0,..,k_N^0)\in(0,1)^N$
such that ${\bf A}_0({\bf k}^0){\bf f}\ne0$ for some constant ${\bf
f}\in\C^M$. Consider the functions
\[\lb{100}
 \eta_{i}(k_i)=\frac1{\sqrt{2\ve}}\ca1,& k_i\in[k_i^0-\ve,k_i^0+\ve],\\ 0,& otherwise,\ac
\]
where $\ve>0$ is some small number. Define $\eta({\bf
k})=\prod_{i=1}^N\eta_i(k_i)$. It is not difficult to check that
$\langle{\bf A}_{\a}{\bf f}\eta\rangle_{\a}$ (where $\a\ne\es$) are
small for small $\ve$ because integrals $\langle\eta\rangle_{\a}$
are small. Due to ${\bf A}_0({\bf k}^0){\bf f}\ne0$ we have $\cA{\bf
f}\eta={\bf A}_{\es}{\bf f}\eta+\sum_{\a\in I_N\sm\es}\langle{\bf
A}_{\a}{\bf f}\eta\rangle_{\a}$ is not small for all small $\ve$.
This is contrary to $\cA=0$. Then ${\bf A}_{\es}=0$. On the next
steps we consistently obtain ${\bf A}_{\{i\}}=0$, after that ${\bf
A}_{\{i,j\}}=0$ and so on. For this we only need to take
$\eta_{\a}=\prod_{i\in\a}\eta_i$ instead of $\eta$ and use the same
arguments as above. \BBox

{\bf Proof of Theorem \ref{T1}.} For any matrix ${\bf
A}=(a_{i,j})_{i,j=1}^M\in\C^{M\ts M}$ denote
\[\lb{101}
 \|{\bf A}\|_{\iy}=\max_{i}\sum_{j=1}^M|a_{i,j}|.
\]
It is a norm satisfying $\|{\bf A}{\bf B}\|_{\iy}\le\|{\bf
A}\|_{\iy}\|{\bf B}\|_{\iy}$. The function \er{004} is
\[\lb{102}
 \|\cA\|_{\mL}=\sum_{\a\in I_N}\|\langle{\bf A}_{\a}({\bf k},{\bf
 x}_{\a})\cdot\rangle_{\a}\|_{\mL}=\sum_{\a\in I_N}\max_{({\bf k},{\bf
 x}_{\a})}\|{\bf A}({\bf k},{\bf x}_{\a})\|_{\iy}.
\]
It is obvious that $\|\l\cA\|_{\mL}=|\l|\|\cA\|_{\mL}$. Due to Lemma
\ref{L1} we have $\|\cA\|_{\mL}=0$ if and only if $\cA=0$. Consider
two operators $\cA,\cB\in\mL$
\[\lb{103}
 \cA=\sum_{\a\in I_N}\langle{\bf A}_{\a}({\bf k},{\bf
 x}_{\a})\cdot\rangle_{\a},\ \ \cB=\sum_{\a\in I_N}\langle{\bf B}_{\a}({\bf k},{\bf
 x}_{\a})\cdot\rangle_{\a}.
\]
Then
\[\lb{104}
 \|\cA+\cB\|_{\mL}=\sum_{\a}\max_{({\bf k},{\bf
 x}_{\a})}\|({\bf A}+{\bf B})({\bf k},{\bf x}_{\a})\|_{\iy}\le\sum_{\a}\max_{({\bf k},{\bf
 x}_{\a})}(\|{\bf A}({\bf k},{\bf x}_{\a})\|_{\iy}+\|{\bf B}({\bf k},{\bf
 x}_{\a})\|_{\iy})
\]
\[\lb{105}
 \le\sum_{\a}(\max_{({\bf k},{\bf
 x}_{\a})}\|{\bf A}({\bf k},{\bf x}_{\a})\|_{\iy}+\max_{({\bf k},{\bf
 x}_{\a})}\|{\bf B}({\bf k},{\bf
 x}_{\a})\|_{\iy})=\|\cA\|_{\mL}+\|\cB\|_{\mL}.
\]
Due to \er{no1}, \er{no2}, \er{102} we have
\[\lb{106}
 \|\langle{\bf A}({\bf k},{\bf x}_{\a})\cdot\rangle_{\a}\langle{\bf B}({\bf k},{\bf
 x}_{\b})\cdot\rangle_{\b}\|_{\mL}=\|\langle({\bf A}\circ{\bf B})({\bf k},{\bf
 x}_{\a\cup\b})\cdot\rangle_{\a\cup\b}\|_{\mL}
\]
\[\lb{107}
 =\max_{({\bf k},{\bf
 x}_{\a\cup\b})}\|\int_{[0,1)^{|\a\cap\b|}}{\bf
 A}({\bf k},{\bf x}_{\a\sm\b}\diamond{\bf z}_{\a\cap\b})
 {\bf B}({\bf k}_{\ol{\a}}\diamond{\bf z}_{\a\cap\b}\diamond{\bf x}_{\a\sm\b},{\bf
 x}_{\b})d{\bf z}_{\a\cap\b}\|_{\iy}
\]
\[\lb{108}
 \le\max_{({\bf k},{\bf
 x}_{\a\cup\b})}\int_{[0,1)^{|\a\cap\b|}}\|{\bf
 A}({\bf k},{\bf x}_{\a\sm\b}\diamond{\bf z}_{\a\cap\b})
 {\bf B}({\bf k}_{\ol{\a}}\diamond{\bf z}_{\a\cap\b}\diamond{\bf x}_{\a\sm\b},{\bf
 x}_{\b})\|_{\iy}d{\bf z}_{\a\cap\b}
\]
\[\lb{109}
 \le\max_{({\bf k},{\bf
 x}_{\a\cup\b})}\max_{{\bf z}_{\a\cap\b}}\|{\bf
 A}({\bf k},{\bf x}_{\a\sm\b}\diamond{\bf z}_{\a\cap\b})
 {\bf B}({\bf k}_{\ol{\a}}\diamond{\bf z}_{\a\cap\b}\diamond{\bf x}_{\a\sm\b},{\bf
 x}_{\b})\|_{\iy}
\]
\[\lb{110}
 \le\max_{({\bf k},{\bf
 x}_{\a\cup\b})}\max_{{\bf z}_{\a\cap\b}}\|{\bf
 A}({\bf k},{\bf x}_{\a\sm\b}\diamond{\bf z}_{\a\cap\b})\|_{\iy}\|
 {\bf B}({\bf k}_{\ol{\a}}\diamond{\bf z}_{\a\cap\b}\diamond{\bf x}_{\a\sm\b},{\bf
 x}_{\b})\|_{\iy}
\]
\[\lb{111}
 \le\max_{({\bf k},{\bf
 x}_{\a})}\|{\bf A}({\bf k},{\bf x}_{\a})\|_{\iy}\max_{({\bf k},{\bf
 x}_{\b})}\|{\bf B}({\bf k},{\bf x}_{\b})\|_{\iy}=\|\langle{\bf A}({\bf k},{\bf x}_{\a})\cdot\rangle_{\a}\|_{\mL}
 \|\langle{\bf B}({\bf k},{\bf
 x}_{\b})\cdot\rangle_{\b}\|_{\mL}
\]
Due to \er{no1} we have
\[\lb{112}
 \cA\cB=\cC=\sum_{\a\in I_N}\cC_{\a},\ \ \cC_{\a}=\sum_{\g,\d:\ \g\cup\d=\a}
 \langle{\bf A}_{\g}\cdot\rangle_{\g}\langle{\bf B}_{\d}\cdot\rangle_{\d}.
\]
Using \er{102}, \er{112}, \er{106}-\er{111} we obtain
\[\lb{113}
 \|\cA\cB\|_{\mL}\le\sum_{\a}\|\cC_{\a}\|_{\mL}\le\sum_{\a}\sum_{\g,\d:\
 \g\cup\d=\a}\|\langle{\bf A}_{\g}\cdot\rangle_{\g}\langle{\bf
 B}_{\d}\cdot\rangle_{\d}\|_{\mL}
\]
\[\lb{114}
 \le\sum_{\a}\sum_{\g,\d:\
 \g\cup\d=\a}\|\langle{\bf A}_{\g}\cdot\rangle_{\g}\|_{\mL}\|\langle{\bf
 B}_{\d}\cdot\rangle_{\d}\|_{\mL}=(\sum_{\a}\|\langle{\bf A}_{\a}\cdot\rangle_{\a}\|_{\mL})
 (\sum_{\a}\|\langle{\bf
 B}_{\a}\cdot\rangle_{\a}\|_{\mL})
\]
\[\lb{115}
 =\|\cA\|_{\mL}\|\cB\|_{\mL}.
\]
The identities \er{113}-\er{115} show that $\mL$ is a Banach
algebra. \BBox

\begin{lemma}\lb{L2}
i) If $\cA={\bf A}_{\es}({\bf k})\cdot\in\mH_{\es}$ is invertible
then ${\bf A}_{\es}({\bf k})$ is invertible $\forall{\bf k}$ and
$\cA^{-1}={\bf A}_{\es}^{-1}({\bf k})\cdot\in\mH_{\es}$. ii) If
$\cA\in(\cI+\mL_{\a})$ is invertible in $\mB$ then
$\cA^{-1}\in(\cI+\mL_{\a})$ as well.
\end{lemma}
{\it Proof.} The statement i) is simple, see, e.g. \cite{Kjmaa}.
While the case of continuous kernels is considered in \cite{Kjmaa},
the proof for the staircase kernels is the same. We just note that
if ${\bf A}_{\es}$ is $h$-staircase then ${\bf A}_{\es}^{-1}({\bf
k})$ is also $h$-staircase.

ii) It is true that $\cA=\cI+\langle{\bf A}({\bf k},{\bf x}_{\a})
\cdot\rangle_{\a}$ for some ${\bf A}$. Because ${\bf A}$ is
$h$-staircase, we can separate variables ${\bf A}({\bf k},{\bf
x}_{\a})={\bf C}({\bf k}){\bf D}({\bf x}_{\a})$ with some
$h$-staircase  ${\bf C}$, ${\bf D}$. Note that the matrices ${\bf
C}$, ${\bf D}$ can be non-square. For the operators with separable
kernels there is an explicit representation for the inverse
operators (see, e.g., \cite{Kjmaa})
\[\lb{L2001}
 \cA^{-1}=(\cI+{\bf C}({\bf k})\langle{\bf D}({\bf
 x}_{\a})\cdot\rangle_{\a})^{-1}=\cI-{\bf C}({\bf k}){\bf E}^{-1}({\bf k}_{\ol\a})\langle{\bf D}({\bf
 x}_{\a})\cdot\rangle_{\a},
\]
where
\[\lb{L2002}
 {\bf E}({\bf k}_{\ol\a})={\bf I}+\langle{\bf D}({\bf
 x}_{\a}){\bf C}({\bf k}_{\ol{\a}}\diamond{\bf x}_{\a})\rangle_{\a}.
\]
Note that for the invertibility it is necessary and sufficient to
have $\det{\bf E}({\bf k}_{\ol{\a}})\ne0$, $\forall{\bf
k}_{\ol{\a}}$, see, e.g., \cite{Kjmaa}. It is obvious that ${\bf
E}({\bf k}_{\ol{\a}})$ and ${\bf E}^{-1}({\bf k}_{\ol{\a}})$ are
both staircase. \BBox

\begin{lemma}\lb{L3}
i) Suppose that $\cA=\sum_{\a\in I_N}\cA_{\a}$ is invertible. Then
$\cA_{\es}$ is invertible. ii) Suppose that $\cA=\cI+\sum_{\a\in
J}\cA_{\a}$ is invertible, where $J\ss I_N\sm\es$ and
$\cA_{\a}\in\mH_{\a}$. Then $\cI+\cA_{\a}$ is invertible for each
$\a\in J$ such that $|\a|=\min_{\b\in J}|\b|$.
\end{lemma}
{\it Proof.} i) The proof repeats the arguments of the proof of
Lemma \ref{L1}. We will use the same notations as in Lemma \ref{L1}
and \ref{L2}. We see that if $\cA_{\es}$ is non-invertible then by
Lemma \ref{L2} the matrix ${\bf A}_{\es}({\bf k}^0)$ is
non-invertible for some ${\bf k}^0$ lying strictly inside the cube
of the volume $h^N$ (recall that all functions are $h$-staircase).
Hence there is ${\bf f}\in\C^M$ such that ${\bf A}_{\es}({\bf
k}^0){\bf f}=0$. Then $\cA{\bf f}\eta$ ($\eta$ is defined in Lemma
\ref{L1}) can be arbitrary small (for small $\ve>0$, see Lemma
\ref{L1}) while $\|{\bf f}\eta\|_{L^2}$ is fixed. This means that
$\cA$ is non-invertible (by Banach's Open Mapping Theorem). For ii)
the proof is the same as for i) but for the zero vector we need to
take the zero vector ${\bf f}$ of ${\bf E}({\bf k}^0_{\ol{\a}})$
multiplied by ${\bf C}({\bf k})$ and $\eta_{\ol{\a}}$ (see notations
in Lemmas \ref{L1} and \ref{L2}). We have that if $\cI+\cA_{\a}$ is
non-invertible then $(\cI+\cA_{\a}){\bf C}({\bf k}){\bf
f}\eta_{\ol{\a}}$ can be arbitrary small (for small $\ve>0$, see
Lemma \ref{L1}) while the norm of ${\bf C}({\bf k}){\bf
f}\eta_{\ol{\a}}$ is not small. The same property will be for
$\cA{\bf C}({\bf k}){\bf f}\eta_{\ol{\a}}$ since the norms of
integrals $\langle\eta_{\ol{\a}}\rangle_{\b}$ are small for
$|\b|\ge|\a|$ and $\b\ne\a$. All these arguments are the same as in
\cite{Kjmaa}, where they are presented with more details. \BBox


\begin{lemma}\lb{L5} Suppose that $\prod_{\a\in I_N}\cG_{\a}=\prod_{\a\in
I_N}\wt\cG_{\a}$, where $\cG_{\a},\wt\cG_{\a}\in\mG_{\a}$ and the
orders of the terms in both products are the same. Then
$\cG_{\a}=\wt\cG_{\a}$ for all $\a$.
\end{lemma}
{\it Proof.} Denote $\cG_{\a}=\cI+\cA_{\a}$,
$\wt\cG_{\a}=\cI+\wt\cA_{\a}$, where
$\cA_{\a},\wt\cA_{\a}\in\mL_{\a}$, $\a\ne\es$. Expanding the
products we obtain
\[\lb{L5001}
 \cG_{\es}+\sum_{\a\in I_N\sm\es}\cB_{\a}=\wt\cG_{\es}+\sum_{\a\in I_N\sm\es}\wt\cB_{\a}
\]
with some $\cB_{\a},\wt\cB_{\a}\in\mL_{\a}$. Using Lemma \ref{L1} we
obtain that $\cG_{\es}=\wt\cG_{\es}$. Moreover, all other terms are
also equal $\cB_{\a}=\wt\cB_{\a}$. We have
\[\lb{L5002}
 \cG_{\es}\cA_{\{i\}}=\cB_{\{i\}}=\wt\cB_{\{i\}}=\wt\cG_{\es}\wt\cA_{\{i\}}\ or\
 \cA_{\{i\}}\cG_{\es}=\cB_{\{i\}}=\wt\cB_{\{i\}}=\wt\cA_{\{i\}}\wt\cG_{\es}
\]
depending on the order of terms in the product. Multiplying
\er{L5002} by $\cG_{\es}^{-1}=\wt\cG_{\es}^{-1}$ leads to
$\cA_{\{i\}}=\wt\cA_{\{i\}}$ and hence $\cG_{\{i\}}=\wt\cG_{\{i\}}$
for all $i=1,...,N$. The same arguments allow us to prove
$\cG_{\{i,j\}}=\wt\cG_{\{i,j\}}$ and so on. \BBox

{\bf Proof of Theorem \ref{T3}.} First identity in \er{013} follows
from \er{no1}, second is simple. The statement that $\mL_{\a}$ are
algebras is also simple, and $\mL_{\a}\cap\mL_{\b}=\{0\}$,
$\mG_{\a}\cap\mG_{\b}=\{\cI\}$ for $\a\ne\b$ follow from Lemma
\ref{L1}. The statement that $\mG_{\a}$ is a group follows from
Lemma \ref{L2}. The closedness of $\mG_{\a}$ and $\mH_{\a}$ is
obvious. Suppose that $\cA=\sum_{\a\in I_N}\cA_{\a}$ is invertible,
where $\cA_{\a}\in\mH_{\a}$. By Lemma \ref{L3}
 $\cA_{\es}$ is invertible. Then $\cA_1=\cA_{\es}^{-1}\cA=\cI+\sum_{\a\in
 I_N\sm\es}\cA_{\a1}$ is invertible, where
 $\cA_{\a1}=\cA_{\es}^{-1}\cA_{\a}\in\mH_{\a}$ (see \er{no2}). By Lemma \ref{L3}
 $\cI+\cA_{\{1\}1}$ is invertible with $(\cI+\cA_{\{1\}1})^{-1}=\cI+\cB_{\{1\}}\in\mG_{\{1\}}$
 (note that $\cA_{\{1\}1}$ can be zero).
 Then $(\cI+\cB_{\{1\}})\cA_1=\cI+\sum_{\a\in
 I_N\sm(\es\cup\{1\})}\cA_{\a2}$ is invertible, where
 $\cA_{\a2}\in\mH_{\a}$. Repeating these steps $2^N$ times we consistently
remove the terms with indices $\es$, $\{i\}$ $(i=1,..,N)$, $\{i,j\}$
$(i,j=1,..,N)$ and so on. Note that the order of the steps is
important in the sense that the term with index $\a$ should be
removed after the term with index $\b$ if $|\a|>|\b|$. If
$|\a|=|\b|$ then the order of removing the terms is not important.
Finally, we obtain that $\cA_{2^N}$ defined by
\[\lb{133}
 \cA_{2^N}=\cI+\cA_{\{1,..,N\}2^N},\ \
 \cA_{\{1,..,N\}2^N}\in\mL_{\{1,..,N\}}
\]
is invertible with
\[\lb{134}
 (\cA_{2^N})^{-1}=\cI+\cB_{\{1,..,N\}},\ \
 \cB_{\{1,..,N\}}\in\mL_{\{1,..,N\}}.
\]
Going through all the steps we obtain that
\[\lb{135}
 \cA^{-1}=\prod_{\a\in I_N\sm\es}(\cI+\cB_{\a})\cA_{\es}^{-1},\ \ \cA=\cA_{\es}\prod_{\a\in
 I_N\sm\es}(\cI+\cA_{\a i_{\a}}),\ \
 (\cI+\cA_{\a i_{\a}})^{-1}=\cI+\cB_{\a},
\]
where the orders of terms in the products are the direct and inverse
orders of the steps.
The uniqueness follows from Lemma \ref{L5}. \BBox

\begin{lemma}\lb{L6}
The mapping $\pmb{\t}:\mL\to\mC$ is linear with
$\|\pmb{\t}\|_{\mL\to\mC}=M$ and
$\pmb{\t}(\cA\cB)=\pmb{\t}(\cB\cA)$.
\end{lemma}
{\it Proof.} The linearity of $\pmb{\t}$ is obvious. Due to the
linearity and "multiplicativity" of the standard trace ${\rm Tr}$ of
square matrices, \er{no1}, \er{no2}, \er{008}, \er{002} and the
facts that
\[\lb{116}
 \ol{\a}\sm(\a\cup\b)=\ol{\a\cup\b},\ \
 \ol{\a}\cap(\a\cup\b)=\b\sm\a,\ \ (\a\sm\b)\cup(\a\cap\b)=\a
\]
we obtain
\[\lb{117}
 \t_{\a\cup\b}(\langle{\bf A}({\bf k},{\bf x}_{\a})\cdot\rangle_{\a}\langle{\bf B}({\bf k},{\bf
 x}_{\b})\cdot\rangle_{\b})=\t_{\a\cup\b}(\langle({\bf A}\circ{\bf B})({\bf k},{\bf
 x}_{\a\cup\b})\cdot\rangle_{\a\cup\b})
\]
\[\lb{118}
 ={\rm Tr}\int\limits_{[0,1)^{|\a\cup\b|}}\int\limits_{[0,1)^{|\a\cap\b|}}{\bf
 A}({\bf k}_{\ol{\a\cup\b}}\diamond{\bf x}_{\a\cup\b},{\bf x}_{\a\sm\b}\diamond{\bf z}_{\a\cap\b})
 {\bf B}({\bf k}_{\ol{\a\cup\b}}\diamond{\bf x}_{\b\sm\a}\diamond{\bf z}_{\a\cap\b}\diamond{\bf x}_{\a\sm\b},{\bf
 x}_{\b})d{\bf z}_{\a\cap\b}d{\bf x}_{\a\cup\b}
\]
\[\lb{119}
 ={\rm Tr}\int\limits_{[0,1)^{|\a\cup\b|}}{\bf
 A}({\bf k}_{\ol{\a\cup\b}}\diamond{\bf x}_{\a\cup\b},{\bf x}_{\a})
 {\bf B}({\bf k}_{\ol{\a\cup\b}}\diamond{\bf x}_{\a\cup\b},{\bf
 x}_{\b})d{\bf x}_{\a\cup\b}
\]
\[\lb{120}
 ={\rm Tr}\int\limits_{[0,1)^{|\a\cup\b|}}
 {\bf B}({\bf k}_{\ol{\a\cup\b}}\diamond{\bf x}_{\a\cup\b},{\bf
 x}_{\b}){\bf A}({\bf k}_{\ol{\a\cup\b}}\diamond{\bf x}_{\a\cup\b},{\bf x}_{\a})d{\bf x}_{\a\cup\b}
\]
\[\lb{121}
 =\t_{\a\cup\b}(\langle{\bf B}({\bf k},{\bf
 x}_{\b})\cdot\rangle_{\b}\langle{\bf A}({\bf k},{\bf
 x}_{\a})\cdot\rangle_{\a}).
\]
Due to the linearity of $\pmb{\t}$, \er{112} and \er{117}-\er{121}
we obtain that $\pmb{\t}(\cA\cB)=\pmb{\t}(\cB\cA)$. Definition of
the norm \ref{D2} and $\|\pmb{\t}(\cI)\|_{\mC}=M$ immediately lead
to $\|{\pmb{\t}}\|_{\mL\to\mC}=M$. \BBox

\begin{lemma}\lb{L7} The mapping $\pmb{\pi}$ is multiplicative.
\end{lemma}
{\it Proof.} Consider the function $\wt{\pmb{\pi}}=(\wt\pi_{\a})$
defined by
\[\lb{010}
 \ln\wt{\pi}_{\a}(\cI+\cA)=-\sum_{n=1}^{+\iy}\frac{(-1)^n\t_{\a}(\cA^n)}n
\]
for $\cA\in\mL$ with $\|\cA\|_{\mL}<1$. The inequalities (see also
Theorem \ref{T1} and Lemma \ref{L6})
\[\lb{122}
 \|\pmb{\t}(\cA^n)\|_{\mC}\le M\|\cA^n\|_{\mL}\le M\|\cA\|_{\mL}^n
\]
show that the series \er{010} converges absolutely. Using \er{010},
${\t}_{\a}(\cA\cB)={\t}_{\a}(\cB\cA)$ (see Lemma \ref{L6}), and
repeating the arguments given for the standard traces and
determinants, see, e.g., \cite{GGK}, we obtain
\[\lb{123}
 \ln\wt\pi_{\a}(\cI+\cA)+\ln\wt\pi_{\a}(\cI+\cB)=-\sum_{n=1}^{+\iy}\frac{(-1)^n\t_{\a}(\cA^n)+(-1)^n\t_{\a}(\cB^n)}{n}
\]
\[\lb{124}
 =-\sum_{n=1}^{+\iy}\frac{(-1)^n\t_{\a}((\cA+\cB+\cA\cB)^n)}{n}=\ln\wt\pi_{\a}((\cI+\cA)(\cI+\cB))
\]
if the norms of $\cA$, $\cB$, $\cA+\cB+\cA\cB$ less than $1$. The
properties \er{123}, \er{124} show us that $\wt{\pmb{\pi}}$ is
multiplicative. Moreover,
$\wt{\pmb{\pi}}(\cI+\cA_{\a})=\pmb{\pi}(\cI+\cA_{\a})$ for small
$\cA_{\a}\in\mL_{\a}$ because $\t_{\a}$ is the standard trace of
finite-rank operators and ${\pi}_{\a}$ is the standard determinant
of finite-rank operators. Now consider $\cI+\l\cA$ for $\l\in\C$.
Going through the procedure from the proof of Theorem \ref{T3} (see
\er{L2001},\er{L2002}, \er{133}-\er{135}) we obtain that
\[\lb{125}
 \cI+\l\cA=\prod_{\a\in I_N}(\cI+\langle{\bf A}_{\a}(\l,{\bf k},{\bf
 x}_{\a})\cdot\rangle_{\a}),
\]
where the function ${\bf A}_{\a}$ depends on $\l\in\C$ as a rational
function (here we use the fact that all functions are $h$-staircase
in variables ${\bf k},{\bf x}$ and hence the integrals are just
finite sums). Moreover ${\bf A}_{\a}(\l,{\bf k},{\bf x}_{\a})=0$ for
$\l=0$ because $\cI+\l\cA=\cI$ for $\l=0$ and the representation it
as the product is unique, see Lemma \ref{L5}. Then for small
$\l\in\C$ we have that
\[\lb{126}
 \wt{\pmb{\pi}}(\cI+\l\cA)=(\wt{\pi}_{\a}(\cI+\langle{\bf A}_{\a}(\l,{\bf k},{\bf
 x}_{\a})\cdot\rangle_{\a}))=({\pi}_{\a}(\cI+\langle{\bf A}_{\a}(\l,{\bf k},{\bf
 x}_{\a})\cdot\rangle_{\a}))=\pmb{\pi}(\cI+\l\cA).
\]
Due to the rational dependence on $\l\in\C$  we obtain that \er{126}
is valid for all $\l\in\C$ for which $\|\l\cA\|_{\mL}<1$. Let $\cA$,
$\cB$ be two invertible operators. Using the same arguments as above
and the fact that $\wt{\pmb{\pi}}$ is multiplicative we obtain that
for small $\l\in\C$
\[\lb{127}
 \pmb{\pi}((\cI+\l(\cA-\cI))(\cI+\l(\cB-\cI)))=\wt{\pmb{\pi}}((\cI+\l(\cA-\cI))(\cI+\l(\cB-\cI)))=
\]
\[\lb{127a}
 \wt{\pmb{\pi}}(\cI+\l(\cA-\cI))\wt{\pmb{\pi}}(\cI+\l(\cB-\cI))=\pmb{\pi}(\cI+\l(\cA-\cI))\pmb{\pi}(\cI+\l(\cB-\cI))
\]
which is
\[\lb{128}
 \pmb{\pi}((\cI+\l(\cA-\cI))(\cI+\l(\cA-\cI)))=\pmb{\pi}(\cI+\l(\cA-\cI))\pmb{\pi}(\cI+\l(\cA-\cI))
\]
for small $\l\in\C$. Due to the rational dependence on $\l\in\C$ we
obtain that \er{128} is valid for all $\l\in\C$ for which
$(\cI+\l(\cA-\cI))$ and $(\cI+\l(\cA-\cI))$ are invertible, e.g. for
$\l=1$. Then $\pmb{\pi}(\cA\cB)=\pmb{\pi}(\cA)\pmb{\pi}(\cB)$. \BBox

{\bf Proof of Theorem \ref{T2}.} Many of the results follow from
Lemmas \ref{L6} and \ref{L7}. The identity
$\pmb{\pi}(e^{\cA})=e^{\pmb{\t}(\cA)}$ follows from \er{010} and
$\pmb{\pi}=\wt{\pmb{\pi}}$. The continuity of $\pmb{\pi}$ in the
vicinity of $\cI$ follows from the expansion \er{010}. If
$\cA\in{\rm Inv}(\mL)$ then the continuity in the vicinity of $\cA$
follows from the continuity in the vicinity of $\cI$ because the
multiplicative property gives us
\[\lb{129}
 \pmb{\pi}(\cA+\cB)-\pmb{\pi}(\cA)=\pmb{\pi}(\cA)((\pmb{\pi}(\cI+\cA^{-1}\cB))-\pmb{\pi}(\cI))\
 and\ \|\cA^{-1}\cB\|_{\mL}\le\|\cA^{-1}\|_{\mL}\|\cB\|_{\mL}
\]
is small for small $\cB$. \BBox

\section*{Acknowledgements}

This work was partially supported by Trr 181 project.


\bibliography{bibl_perp1}

\end{document}